
\tolerance=10000
\raggedbottom

\baselineskip=15pt
\parskip=1\jot

\def\sk{\vskip 3\jot}

\def\heading#1{\vskip3\jot{\noindent\bf #1}}
\def\label#1{{\noindent\it #1}}
\def\QED{\hbox{\rlap{$\sqcap$}$\sqcup$}}


\def\ref#1;#2;#3;#4;#5.{\item{[#1]} #2,#3,{\it #4},#5.}
\def\refinbook#1;#2;#3;#4;#5;#6.{\item{[#1]} #2, #3, #4, {\it #5},#6.} 
\def\refbook#1;#2;#3;#4.{\item{[#1]} #2,{\it #3},#4.}


\def\({\bigl(}
\def\){\bigr)}


\def\de{\delta}
\def\ep{\varepsilon}

\def\si{\sigma}
\def\ta{\tau}



\baselineskip=24pt

\def\abs#1{\vert#1\vert}
\def \bigabs#1{\bigg\vert#1\bigg\vert}
\def\norm#1{\Vert#1\Vert}
\def \bignorm#1{\bigg\Vert#1\bigg\Vert}

\def\[{\big[}
\def\]{\big]}

\def\Ex{{\rm Ex}}
\def\Var{{\rm Var}}

\def\wlln{Weak Law of Large Numbers}
\def\clt{Central Limit Theorem}

\def\abs#1{\big\vert#1\big\vert}
\def\norm#1{\vert#1\vert}

{
\pageno=0
\nopagenumbers
\rightline{\tt clt.arxiv.tex}
\vskip1in

\centerline{\bf Elementary Proofs of the Main Limit Theorems of Probability}
\vskip0.5in

\centerline{Nicholas Pippenger}
\centerline{\tt njp@math.hmc.edu}
\sk

\centerline{Department of Mathematics}
\centerline{Harvey Mudd College}
\centerline{1250 Dartmouth Avenue}
\centerline{Claremont, CA 91711}
\vskip0.5in

\noindent{\bf Abstract:}
We give simple proofs, under minimal hypotheses, of the Weak Law of Large Numbers and the Central Limit Theorem for independent identically distributed random variables.
These proofs use only the elementary calculus, together with the most basic notions of probability, expectations, and distribution functions.
\sk
\sk
\vfill\eject
}

\heading{1. Introduction}

Our goal in this paper is to give  simple proofs of the \wlln{} and the \clt{} for independent and identically distributed random variables.
Our proofs make only the minimal necessary assumptions: that the mean exists for the \wlln, and that the mean and variance exist for the \clt.
We use only straightforward manipulations of probabilities, expectations and distribution functions.
In particular, we make no use of characteristic functions, or of other transform techniques, nor does it use any operator-theoretic methods.
We use the Stieltjes integral $\Ex[X] = \int x\,dF_X(x)$ (where $F_X(x) = \Pr[X\le x]$ is the distribution function of $X$) as the definition of expectation.
Our proof holds with complete generality if this integral is interpreted as a Lebesgue-Stieltjes integral.
But everything we do will also be valid if one interprets this definition as the combination of a Riemann integral (for the continuous component of the distribution) and a sum (for the discrete component), ignoring the possibility of a singular component.
The only property of the normal distribution that we use is that a sum of independent normally distributed random variables is again normally distributed, with the means and variances added, as is shown in Ross [3, pp.~256].
Finally, the only analysis that we use is the definition of a limit and Taylor's Theorem (with the Lagrange form of the remainder term).

Our proof of the \clt{} is inspired by that of Trotter [2]; in principle, all we have done is to transform those parts of his proof involving operators and other notions from functional analysis into straightforward manipulations of probabilities, expectations and distribution functions.
And Trotter himself says:
``Our proof is in principle the same as that used by Lindeberg [1].''
So our simple proof actually has a long lineage.
The observation that the same technique also yields the \wlln{} appears to be new.

It is in fact the idea of proving both the \clt{} and the \wlln{}  with a single argument that yields the key to the proof.
In proving any convergence result, it is always tempting to use metric, such as the supremum norm,
since then one can use the triangle inequality and other tools from the theory of metric spaces.
But the supremum norm governs uniform convergence, and while uniform convergence does indeed take place in the \clt, there is in general only pointwise convergence in the \wlln, and
pointwise convergence is not governed by any metric.
This suggests the idea of smoothing the distribution functions of the random variables by convolving them with a smooth distribution function (that is, by adding an independent smoothly distributed random variable).
In the case at hand, the random variables we choose to add are beta-distributed with various parameters.

We prove the \clt{} in the following form.

\label{Theorem 1:}
Let $Z, Z_1, \ldots, Z_n$ be independent and identically distributed random variables with
$\Ex[Z] = 0$ and $\Var[Z] = \Ex[Z^2] = 1$.
Let $U_n = Z_1 + \cdots + Z_n$ and let $S_n = U_n/n^{1/2}$.
Let
$N$ be  normally distributed with $\Ex[N] = 0$ and 
$\Var[N] = 1$.
Then $\lim_{n\to\infty} F_{S_n}(t) = F_N(t)$ for all real $t$.

The \clt{} is often stated in a more general form in which the mean and variance of  $Z$ are not assumed to be $0$ and $1$, respectively.
But this form actually reduces to the special case of Theorem 1:
if $\Ex[Z']=\mu$ and $\Var[Z']=\si^2$, applying Theorem 1 with $Z = (Z'-\mu)/\si$ yields the conclusion that the limiting distribution of $(Z'_1 + \cdots + Z'_n)/n^{1/2}$ is normal, with mean 
$\mu$ and variance $\si^2$.

We shall prove the \wlln{} in  the following form.

\label{Theorem 2:}
Let $Z, Z_1, \ldots, Z_n$ be independent and identically distributed random variables with
$\Ex[Z] = 0$ and $\Ex[\abs{Z}] = 1$.
Let $U_n = Z_1 + \cdots + Z_n$ and let $S_n = U_n/n$.
Let
$D$ be  a deterministic random variable with 
$\Pr[D=0]=1$.
Then $\lim_{n\to\infty} F_{S_n}(t) = F_D(t)$ for all real $t$ at which $F_D(t)$ is continuous
(that is, all $t\not=0$).

The \wlln{} is often stated in terms of ``convergence in probability''
(that is, $\Pr[\abs{S_n/n} > \ep]\to 0$ as $n\to 0$ for all $\ep>0$), rather than in terms of ``convergence in  distribution'', as we have done.
But, although convergence in probability is in general stronger than convergence in  distribution,
in the case of convergence to a deterministic value they are equivalent, since
$$\eqalign{
\Pr[\abs{S_n/n} > \ep] &= \Pr[S_n/n > \ep] + \Pr[S_n/n < -\ep] \cr
&\le 1 - F_{S_n/n}(\ep) + F_{S_n/n}(-\ep), \cr
}$$ 
and convergence in distribution yields $F_{S_n/n}(\ep)\to 1$ and $F_{S_n/n}(-\ep)\to 0$ as $n\to \infty$ for all $\ep>0$.

The \wlln{} is often stated in a more general form in which the means  of  $Z$ and $\abs{Z}$ are not assumed to be $0$ and $1$, respectively.
But this form actually reduces to the special case of Theorem 2:
if $\Ex[Z']=\mu$ and $\Ex[\abs{Z'-\mu}] = \ta$, applying Theorem 2 with $Z = (Z-\mu)/\ta$ yields the conclusion that the limiting distribution of $(Z'_1 + \cdots + Z'_n)/n$ is deterministic, with almost sure value $\mu$.
\sk

\heading{2. The Proofs}

Our key tactic, which will be used repeatedly in our proofs, is that if $W$ and $X$ are independent random variables,
$$F_{W+X}(w) = \int F_W(w-x)\,dF_X(x) = \Ex[F_W(w-X)].$$
The strategy of our proof can be described in rough terms as follows.
To prove that a sequence $F_{S_n}$ of distributions tends to a limiting distribution $F_T$, 
it is natural to try to use a metric in the space of distributions, such as
$$\bignorm{F_{S_n} - F_T} = \sup_t \bigabs{F_{S_n}(t) - F_T(t)}.$$
But this metric describes uniform convergence of the distributions, which may not hold in the cases we consider.
(For example,  we will not usually even have convergence at $t=0$ in Theorem 2.)
This leads us to employ the following device.
We add to both $S_n$ and $T$ an independent random variable $W$.
By giving $W$ a sufficiently smooth distribution function, we will be able to prove
that $F_{W+S_n}$ converges uniformly to $F_{W+T}$, and by making $W$ sufficiently small
in absolute value, even pointwise convergence will imply the convergence of $F_{S_n}(t)$ to  
$F_T(t)$
for all $t$ at which $F_T$ is continuous.

\label{Lemma 3:}
Let $S_n$ for $n\ge 1$ and $T$ be random variables.
If, for every $\de>0$, there exists a random variable $W$, independent of the $S_n$ and $T$,
satisfying $\abs{W} \le \de$ and such that
$F_{S_n + W}(s) \to F_{T + W}(s)$ as $n\to \infty$ for all $s$,
then $F_{S_n}(t) \to F_{T}(t)$ as $n\to \infty$ for all $t$ at which $F_T$ is continuous.

\label{Proof:}
Let $t$ be a point of continuity for $F_T$.
Given $\ep > 0$, find $\de > 0$ such that $F_T(t+2\de)\le F_T(t) + \ep$
and $F_T(t-2\de) \ge F_T(t)-\ep$.
Let $W$ satisfy the hypotheses of the lemma.
Since $W\le\de$, we have 
$$\eqalign{
F_{S_n}(t) 
&= \Pr[S_n\le t] \cr
&\le \Pr[S_n + W\le t + \de] \cr
&= F_{S_n + W}(t + \de). \cr
}$$
Since $W\ge-\de$, we have
$$\eqalign{
F_{T+W}(t+\de) 
&= \Pr[T+W\le t+\de] \cr
&= \Pr[T\le t+\de-W] \cr
&\le \Pr[T\le t+2\de] \cr
&= F_T(t+2\de). \cr
}$$
Since $\lim_{n\to\infty} F_{S_n+W}(t+\de) = F_{T+W}(t+\de)$,  we have
$$\eqalign{
F_{S_n(t)}
&\le F_{S_n+W}(t+\de) \cr
&\le F_{T+W}(t+\de) + \ep \cr
&\le F_T(t+2\de) + \ep \cr
&\le F_T(t) + 2\ep \cr
}$$
for all sufficiently large $n$.
A similar argument yields
$$F_{S_n(t)} \ge F_T(t) - 2\ep$$
for all sufficiently large $n$.
Since these inequalities hold for all $\ep>0$, we obtain
$\lim_{n\to\infty} F_{S_n}(t) = F_T(t)$,
as desired.
\QED

If $S_n$ is the sum of $n$ independent and identically distributed contributions $X_1, \ldots, X_n$,
and $T_n$ is the sum of $n$ independent and identically distributed contributions 
$Y_1, \ldots, Y_n$, then we can go from $S_n$ to $T_n$ in $n$ steps by changing one $X_i$ at each step to $Y_i$.
Then $n$ applications of the triangle inequality  will allow us to bound $\norm{F_{S_n} - F_T}$ by $n$ times $\norm{F_X - F_Y}$, where $X$ and $Y$ have the common distributions of $X_i$ and $Y_i$, respectively.
Adding the contribution of $W$ in each case yields the following lemma.

\label{Lemma 4:}
Let $X, X_1, \ldots, X_n$ be independent and identically distributed random variables,
and let $S_n = X_1 + \cdots + X_n$.
Let $Y, Y_1, \ldots, Y_n$ be  identically distributed random variables, 
independent of each other and the $X_i$,
and let $T_n = Y_1 + \cdots + Y_n$.
Let $W$ be a random variable independent of $X$, the $X_i$, $Y$, and the $Y_i$.
Then
$$\bignorm{F_{W+S_n} -  F_{W+T_n}}\le n\,\bignorm{F_{W+X} -  F_{W+Y}}.$$

\label{Proof:}
For $0\le i\le n$, let $Q_i = Y_1 + \cdots + Y_i$ and $R_i = X_i  + \cdots + X_n$.
Then
$$\eqalign{
F_{W+S_n}(w) - &F_{W+T_n}(w) \cr
&= \sum_{1\le i\le n} F_{W+Q_{i-1}+X_i+R_{i+1}}(w) - F_{W+Q_{i-1}+Y_i+R_{i+1}}(w) \cr
&= \sum_{1\le i\le n} \Ex[F_{W+X_i}(w-Q_{i-1}-R_{i+1}) - F_{W+Y_i}(w-Q_{i-1}-R_{i+1})]. \cr
}$$
Taking absolute values, we obtain
$$\eqalign{
\bigabs{F_{W+S_n}(w) - &F_{W+T_n}(w)} \cr
&= \left\vert\sum_{1\le i\le n} 
\Ex[F_{W+X_i}(w - Q_{i-1} - R_{i+1}) - F_{W+Y_i}(w - Q_{i-1} - R_{i+1})] \right\vert \cr
&\le\sum_{1\le i\le n}  
\bigabs{\Ex[F_{W+X_i}(w - Q_{i-1} - R_{i+1}) - F_{W+Y_i}(w - Q_{i-1} - R_{i+1})]} \cr
&\le\sum_{1\le i\le n} 
\Ex\left[\,  \bigabs{F_{W+X_i}(w - Q_{i-1} - R_{i+1}) - F_{W+Y_i}(w - Q_{i-1} - R_{i+1})}\,\right]. \cr
}$$
Each absolute value in the last right-hand side is bounded by $\norm{F_{W+X} - F_{W+Y}}$.
Since there are $n$ terms in the sum, we obtain
$$\bigabs{ F_{W+S_n}(w) - F_{W+T_n}(w)} \le n\, \bignorm{F_{W+X} - F_{W+Y}}.$$
Since this inequality holds for all $w$, we obtain the conclusion of the lemma.
\QED

In Lemma 3, we required $W$ to be small compared with $S_n$ and $T$.
But if $S_n$ and $T_n = T$ are each the sum of $n$ identically distributed contributions $X_i$ and $Y_i$, as in Lemma 4, then these contributions will be small compared to $W$ when $n$ is large.
Thus, if the distribution of $W$ is smooth, it will change little when the small random variable
$X_i$ or $Y_i$ is added.
The following lemma and its corollary express this fact in the form we need.

\label{Lemma 5:}
Let $W$ and $Z$ be independent random variables, with  $F_W$ having two bounded and uniformly continuous derivatives, and with $Z$ satisfying $\Ex[Z]=0$ and $\Ex[Z^2]=1$.
Let $X = Z/N^{1/2}$.
Then, for every $\ep>0$,
$$\bigabs{F_{W+X}(w) - F_W(w) - {1\over 2n}F''_W(w)} \le {\ep\over n}$$
for  all sufficiently large $n$ and all $w$.

\label{Proof:}
We have
$$F_{W+X}(w) = \Ex[F_W(w - X)].$$
Since $F_W$ has two continuous derivatives, we may expand it in a Taylor series,
$$\eqalignno{
F_W(w - X)
&= F_W(w - Z/n^{1/2}) \cr
&= F_W(w)  - {1\over n^{1/2}} F'_W(w)\,Z + {1\over 2n} F''_W(v(Z)) Z^2 \cr
&= F_W(w)  - {1\over n^{1/2}} F'_W(w)\,Z
+ {1\over 2n} F''_W(w) + {1\over 2n} \(F''_W(v(Z)) - F''_W(w)\) Z^2, \cr
}$$
where $v$ is a function satisfying $w - Z/n^{1/2} \le v(Z)\le w$, 
with $F''_W(v(Z)) Z^2$ integrable because all the other
terms in the equation defining it are integrable.
Taking expectations on both sides, and using $\Ex[Z]=0$ and $\Ex[Z^2]=1$, we have
$$F_{W+X}(w) = F_W(w) + {1\over 2n}F''_W(w) 
+ {1\over 2n}\Ex\[ \(F''_W(v(Z)) - F''_W(w)\) Z^2 \].$$
Thus, given $\ep>0$,
it will suffice to show that the last term on the right-hand side has absolute value at most $2\ep$ for all sufficiently large $n$.
Since $F_W$ has a uniformly continuous second derivative, there exists $\de>0$ such that
$\abs{v - w} \le \de$ implies $\abs{F''_W(v) - F''_W(w)} \le \ep$ for all $v$ and $w$.
And since the second derivative is bounded, there exists $M$ such that $\abs{F''_W(v)} \le M$ for all 
$v$.
We break the expectation to be bounded into two parts,
$$\eqalign{
\Ex\[ \(F''_W(v(Z)) - F''_W(w)\) Z^2 \] 
&= \Ex\[ \(F''_W(v(Z)) - F''_W(w)\) Z^2, \abs{Z} \le \de n^{1/2} \] \cr
&\qquad+  \Ex\[ \(F''_W(v(Z)) - F''_W(w)\) Z^2, \abs{Z} > \de n^{1/2} \], \cr
}$$
where $\Ex[g(Z), E] = \int_E g(Z)\,dF_Z(x)$ denotes the expectation of $g(Z)$ restricted to the event
$E$.
Since $\abs{Z} \le \de n^{1/2}$ implies $\abs{v(Z) - w} \le \de$, which in turn implies 
$\abs{F''_W(v(Z)) - F''_W(w)} \le \ep$, we have
$$\eqalign{
\bigabs{\Ex\[ \(F''_W(v(Z)) - F''_W(w)\) Z^2, \abs{Z} \le \de n^{1/2} \] }
&\le \ep\, \Ex[  Z^2, \abs{Z} \le \de n^{1/2} ] \cr
&\le \ep\, \Ex[Z^2] \cr
& \le \ep. \cr
}$$
And since $\abs{F''_W(v(Z)) - F''_W(w)} \le 2M$ and
$\Ex[Z^2, \abs{Z} > y] \to 0$ as $y\to\infty$, we have
$$\eqalign{
\bigabs{\Ex\[ \(F''_W(v(Z)) - F''_W(w)\) Z^2, \abs{Z} > \de n^{1/2} \] }
&\le 2M\, \Ex[  Z^2, \abs{Z} > \de n^{1/2} ] \cr
&\le \ep \cr
}$$ for all sufficiently large $n$.
\QED

\label{Corollary 6:}
Let $W$, $Z$ and $Z'$ be independent random variables, with  $F_W$ having two bounded and uniformly continuous derivatives, and with $Z$ and $Z'$ satisfying $\Ex[Z]=\Ex[Z']=0$ and 
$\Ex[Z^2]=\Ex[Z^{\prime 2}]=1$.
Let $X = Z/n^{1/2}$ and $X' = Z'/n^{1/2}$.
Then, for every $\ep>0$,
$$\bignorm{F_{W+X} - F_{W+X'}(w)} \le {\ep\over n}.$$

\label{Proof:}
Apply Lemma 5 with $\ep/2$ for $\ep$, and $Z$ and $Z'$ in turn for $Z$, and use the triangle inequality for absolute value.
\QED

\label{Proof of Theorem 1:}
Since $F_N$ is continuous, we may use Lemma 3 to show convergence at all points.
Given $\de>0$,
we take $W$ to be the median of five independent random variables, each uniformly distributed
in the interval $[-\de,\de]$ (equivalently, $W = \de(2B - 1)$, where $B\sim {\rm Beta}(3,3)$).
This random variable clearly meets the conditions of Lemma 3 and, since its first two derivative vanish at $\pm\de$ and it varies only over a closed and bounded interval, it also meets the conditions of  Lemma 5 and Corollary 6.
Taking $Z'$ to be the normal random variable $N$ and $Y = N/n^{1/2}$ in Corollary 6, we conclude that, for every 
$\ep>0$,
$$\bignorm{F_{W+X} - F_{W+Y}(w)} \le {\ep\over n}.$$
Lemma 4 then allows us to conclude that,
for every $\ep>0$,
$$\bignorm{F_{W+S_n} -  F_{W+T_n}}\le \ep,$$
where $S_n = X_1 + \cdots + X_n$ (with the $X_i$ independently distributed like $X$) and
$T_n = Y_1 + \cdots + Y_n$ (with the $Y_i$ independently distributed like $Y$)
Since $T_n$ is the sum of $n$ independent and normally distributed random variables,
each having variance $1/n$, it has the distribution of the standard normal random variable $N$.
Thus we obtain
$$\bignorm{F_{W+S_n} -  F_{W+N}}\le \ep.$$
We can now apply Lemma 3 with  $N$ for $T$, and conclude that
$F_{S_n}(t)$ converges to $F_N(t)$ for all $t$.
\QED

Our proof of Theorem 2 is even simpler than our proof of Theorem 1.

\label{Lemma 7:}
Let $W$ and $Z$ be independent random variables, with  $F_W$ having a bounded and uniformly continuous derivative, and with $Z$ satisfying $\Ex[Z]=0$ and $\Ex[\abs{Z}]\le 1$.
Let $X = Z/n$.
Then, for every $\ep>0$,
$$\bigabs{F_{W+X}(w) - F_W(w)} \le {\ep\over n}$$
for  all sufficiently large $n$ and all $w$.

\label{Proof:}
We have
$$F_{W+X}(w) = \Ex[F_W(w - X)].$$
Since $F_W$ has a continuous derivative, we may expand it in a Taylor series,
$$\eqalignno{
F_W(w - X)
&= F_W(w - Z/n) \cr
&= F_W(w)  - {1\over n} F'_W(v(Z)) Z \cr
&= F_W(w) 
+ {1\over n} F'_W(w) Z + {1\over n} \(F'_W(v(Z)) - F'_W(w)\) Z, \cr
}$$
where $v$ is a function satisfying $w - Z/n \le v(Z)\le w$, 
with $F'_W(v(Z)) Z$ integrable because all the other
terms in the equation defining it are integrable.
Taking expectations on both sides, and using $\Ex[Z]=0$, we have
$$F_{W+X}(w) = F_W(w)  + {1\over n}\Ex\[ \(F'_W(v(Z)) - F'_W(w)\) Z \].$$
Thus, given $\ep>0$,
it will suffice to show that the last term on the right-hand side has absolute value at most $2\ep$ for all sufficiently large $n$.
Since $F_W$ has a uniformly continuous derivative, there exists $\de>0$ such that
$\abs{v - w} \le \de$ implies $\abs{F'_W(v) - F'_W(w)} \le \ep$ for all $v$ and $w$.
And since the derivative is bounded, there exists $M$ such that $\abs{F'_W(v)} \le M$ for all 
$v$.
We break the expectation to be bounded into two parts,
$$\eqalign{
\Ex\[ \(F'_W(v(Z)) - F'_W(w)\) Z \] 
&= \Ex\[ \(F'_W(v(Z)) - F'_W(w)\) Z, \abs{Z} \le \de n \] \cr
&\qquad+  \Ex\[ \(F'_W(v(Z)) - F'_W(w)\) Z, \abs{Z} > \de n \]. \cr
}$$
Since $\abs{Z} \le \de n$ implies $\abs{v(Z) - w} \le \de$, which in turn implies 
$\abs{F'_W(v(Z)) - F'_W(w)} \le \ep$, we have
$$\eqalign{
\bigabs{\Ex\[ \(F'_W(v(Z)) - F'_W(w)\) Z, \abs{Z} \le \de n \] }
&\le \ep\, \Ex[  \abs{Z}, \abs{Z} \le \de n ] \cr
&\le \ep\, \Ex[\abs{Z}] \cr
& \le \ep. \cr
}$$
And since $\abs{F'_W(v(Z)) - F'_W(w)} \le 2M$ and
$\Ex[\abs{Z}, \abs{Z} > y] \to 0$ as $y\to\infty$, we have
$$\eqalign{
\bigabs{\Ex\[ \(F'_W(v(Z)) - F'_W(w)\) Z, \abs{Z} > \de n \] }
&\le 2M\, \Ex[  \abs{Z}, \abs{Z} > \de n ] \cr
&\le \ep \cr
}$$ for all sufficiently large $n$.
\QED

\label{Corollary 8:}
Let $W$, $Z$ and $Z'$ be independent random variables, with  $F_W$ having a bounded and uniformly continuous derivative, and with $Z$ and $Z'$ satisfying $\Ex[Z]=\Ex[Z']=0$ and 
$\Ex[\abs{Z}], Ex[\abs{Z'}]\le 1$.
Let $X = Z/n$ and $X' = Z'/n$.
Then, for every $\ep>0$,
$$\bignorm{F_{W+X} - F_{W+X'}(w)} \le {\ep\over n}.$$

\label{Proof:}
Apply Lemma 7 with $\ep/2$ for $\ep$, and $Z$ and $Z'$ in turn for $Z$, and use the triangle inequality for absolute value.
\QED

\label{Proof of Theorem 2:}
Since $F_D$ is continuous at all points other than $0$, we may use Lemma 3 to show convergence at all points other than $0$.
Given $\de>0$,
we take $W$ to be the median of three independent random variables, each uniformly distributed
in the interval $[-\de,\de]$ (equivalently, $W = \de(2B - 1)$, where $B\sim {\rm Beta}(2,2)$).
This random variable clearly meets the conditions of Lemma 3 and, since its first derivative vanishes at $\pm\de$ and it varies only over a closed and bounded interval, it also meets the conditions of  Lemma 7 and Corollary 8.
Taking $Z'$ to be the deterministic random variable $D$ in Corollary 8, and $Y = D/n+ D$,
we conclude that, for every 
$\ep>0$,
$$\bignorm{F_{W+X} - F_{W+Y}(w)} \le {\ep\over n}.$$
Lemma 4 then allows us to conclude that,
for every $\ep>0$,
$$\bignorm{F_{W+S_n} -  F_{W+T_n}}\le \ep.$$
Since $T_n$ is the sum of $n$ deterministic  random variables,
each having mean $0$, it has the distribution of the deterministic random variable $D$.
Thus we obtain
$$\bignorm{F_{W+S_n} -  F_{W+D}}\le \ep.$$
We can now apply Lemma 3 with $D$ for $T$, and conclude that
$F_{S_n}(t)$ converges to $F_D(t)$ for all $t$ at which $F_D$ is continuous
(that is, all $t\not=0$).
\QED
\sk

\heading{3. Acknowledgment.} This research was partially supported by NSF Grant CCF CCF 0917026.
\sk

\heading{4. References}

\ref 1; J. W. Lindeberg;
``Eine neue Herleitung des Exponentialgesetzes in der Wahrscheinlichkeitsrechnung'';
Math.\ Z.; 15 (1922) 211--225.

\ref 2; H. F. Trotter;
``An Elementary Proof of the Central Limit Theorem'';
Arch.\ Math.; 10 (1959) 226--234.

\bye